\documentstyle[fullpage,11pt,amssymb]{article}

\input psfig
\input xypic
\def\SBIMSMark#1#2#3{
 \font\SBF=cmss10 at 10 true pt
 \font\SBI=cmssi10 at 10 true pt
 \setbox0=\hbox{\SBF Stony Brook IMS Preprint \##1}
 \setbox2=\hbox to \wd0{\hfil \SBI #2}
 \setbox4=\hbox to \wd0{\hfil \SBI #3}
 \setbox6=\hbox to \wd0{\hss
             \vbox{\hsize=\wd0 \parskip=0pt \baselineskip=10 true pt
                   \copy0 \break%
                   \copy2 \break%
                   \copy4 \break}}
 \dimen0=\ht6   \advance\dimen0 by \vsize \advance\dimen0 by 8 true pt
                \advance\dimen0 by -\pagetotal
 \dimen2=\hsize \advance\dimen2 by .25 true in
  \openin2=publishd.tex
  \ifeof2\setbox0=\hbox to 0pt{}
  \else 
     \setbox0=\hbox to 3.1 true in{
                \vbox to \ht6{\hsize=3 true in \parskip=0pt  \noindent  
                \input publishd.tex 
                \vfill}}
  \fi
  \closein2
  \ht0=0pt \dp0=0pt
 \ht6=0pt \dp6=0pt
 \setbox8=\vbox to \dimen0{\vfill \hbox to \dimen2{\copy0 \hss \copy6}}
 \ht8=0pt \dp8=0pt \wd8=0pt
 \copy8
 \message{*** Stony Brook IMS Preprint #1, #2 ***}
}

\newtheorem{thm}{Theorem}[section]
\newtheorem{theorem}[thm]{Theorem}
\newtheorem{proposition}[thm]{Proposition}

\newtheorem{lemma}[thm]{Lemma}

\newenvironment{MT}[1]{\bigskip 
\noindent {\bf Theorem} {\bf #1 \/} \begin{em}}{\end{em} \bigskip}

\newenvironment{remark}{\medskip 
\noindent {\bf Remark.}}{\mbox{}\medskip}
 
\newenvironment{definition}{\medskip 
\noindent {\bf Definition.}}{\mbox{}\medskip} 
\newenvironment{proof}{\medskip
\noindent {\it Proof:}}{\hfill $\odot$ \mbox{}\bigskip}

\newcommand{\mod}{\text{mod}}
\newcommand{\dist}{\text{dist}}
\newcommand{\capp}{\text{cap}}

\begin{document}


\begin{center}
{\Large \bf Parameter Scaling for the Fibonacci Point}\\
\vspace{2 cm}
LeRoy Wenstrom\\
Mathematics  Department\\
S.U.N.Y.  Stony Brook\\
\end{center}

\begin{abstract}
We prove geometric and scaling results for the real Fibonacci parameter value in the quadratic family $z \mapsto z^2+c$.  The principal nest of the Yoccoz parapuzzle pieces has rescaled asymptotic geometry equal to the filled-in Julia set of $z \mapsto z^2-1$.  The modulus of two such successive parapuzzle pieces increases at a linear rate.  Finally, we prove a ``hairiness" theorem for the Mandelbrot set at the Fibonacci point when rescaling at this rate.
\end{abstract}

\SBIMSMark{1996/4}{June 1996}{}

In this paper, we focus on the small scale similarities between the dynamical space 
and parameter space for the Fibonacci point in the family of maps $z \mapsto z^2+c$.  There 
is a general philosophy in complex dynamics that the structure we see in the parameter space 
around the parameter value $c$ should be the ``same" as that around the critical value `$c$' 
in dynamical space \cite{DH3}.  In the case where the critical point is pre-periodic, 
Tan Lei \cite{T} proved such asymptotic similarities by showing that the Mandelbrot set 
and Julia set exhibit the same limiting geometry.   For parameters in which the critical 
point is recurrent (i.e., it eventually returns back to any neighborhood of itself), the 
Mandelbrot and Julia sets are much more complicated.  Milnor, in \cite{Mp}, made a number 
of conjectures (as well as pictures!) for the case of infinitely renormalizable
points of bounded type.  Dilating by factors determined by the renormalization, the 
resulting computer pictures demonstrate a kind of self-similarity, with each successive 
picture looking like a ``hairier" copy of the previous.  McMullen \cite{McM2} has proven that, 
for these points, the Julia set densely fills the plane upon repeated rescaling, i.e., hairiness;
and Lyubich has recently proven hairiness of the Mandelbrot set for Feigenbaum like points.  We focus on a primary example 
of dynamics in which we have a recurrent critical point and the dynamics is 
non-renormalizable: the Fibonacci map.  

The dynamics of the real quadratic Fibonacci map, where the critical point returns 
closest to itself at the Fibonacci iterates, has been extensively studied 
(especially see \cite{LM}).  Maps with Fibonacci type returns were first discovered 
in the cubic case by Branner and Hubbard \cite {BH} and have since been consistently 
explored because they are a fundamental combinatorial type of the class of 
non-renormalizable maps. 
The Fibonacci map was used by Lyubich and Milnor in developing the 
{\it generalized renormalization} procedure which has proven very fruitful.  
The Fibonacci map was also highlighted in the work of Yoccoz as it was in some 
sense the worst case in the proof of local connectivity of non-renormalizable 
Julia sets with recurrent critical point \cite{H}, \cite{Ml}. 

The local connectedness proof of Yoccoz involves producing a sequence of partitions 
of the Julia set, now called Yoccoz puzzle pieces.  These Yoccoz puzzle pieces are 
then shown to exhibit the {\it divergence property} and in particular nest down to 
the critical point, proving local connectivity there.  Yoccoz then transfers this 
divergence property to the parapuzzle pieces around the parameter point to demonstrate 
that the Mandelbrot set is locally connected at this parameter value.   Lyubich further
explores the Yoccoz puzzle pieces of Fibonacci maps and 
demonstrates that the {\it principal nest} of Yoccoz puzzle pieces has rescaled 
asymptotic geometry equal to the filled-in Julia set of $z \mapsto z^2-1$ and that 
the moduli of successive annuli grow at a linear rate \cite{Lt}.

We prove that the same geometric and rescaling results hold for the principal nest 
of parapuzzle pieces for the Fibonacci parameter point in the Mandelbrot set.  Let 
the notation $mod(A,B)$ (where $B \subset A$) indicate the modulus of the annulus 
$A \setminus B$. (See Appendix for the definition of the modulus.)

\begin{MT}{A:}   {\bf (Parapuzzle scaling and geometry)}     \label{MT}

The principal nest of Yoccoz parapuzzle pieces, $P^n$, for the Fibonacci point $c_{fib}$ 
has the following properties. \\
1. They scale down to the point $c_{fib}$ in the following asymptotic manner:

\[  \lim_{n \rightarrow \infty} \mod (P^{n-1}, P^{n}) \;  / \;  n  =
\frac{2}{3}\ln 2.   \]
2.  The rescaled $P^n$ and the boundary of the rescaled $P^n$ have
asymptotic geometry equal to
the
filled-in Julia set of $z \mapsto z^2-1$ and its boundary, respectively.   
\end{MT}
By definition, the sets $P^n$, suitably rescaled, have asymptotic geometry
equal to a set $K$ if there are complex affine transformations $A_n$ so
that the
images $A_n(P^n)$ converge to $K$ in the Hausdorff metric.
 
\begin{figure}[hbt]       
\centerline{\psfig{figure=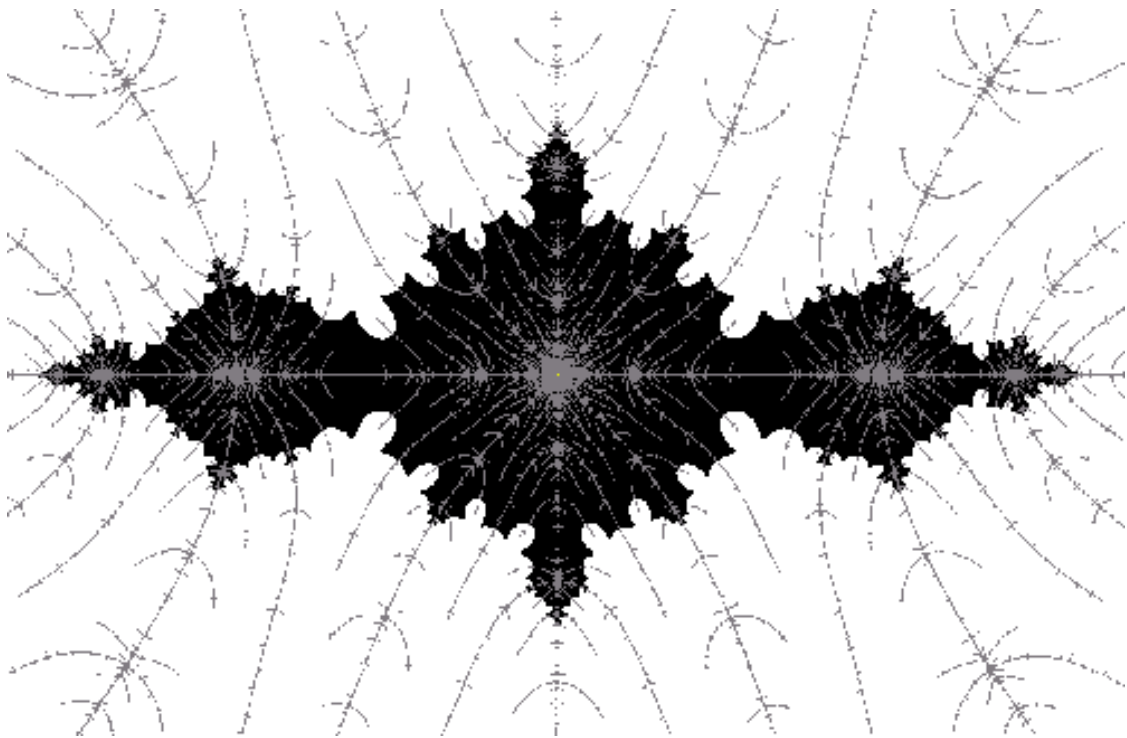,width=.9\hsize}}
\caption{$P^7$ for the Fibonacci point, the seventh level parapuzzle piece with a part of the Mandelbrot set. }
\label{Param puzzle}
\end{figure}

\begin{figure}[hbt]            
\centerline{\psfig{figure=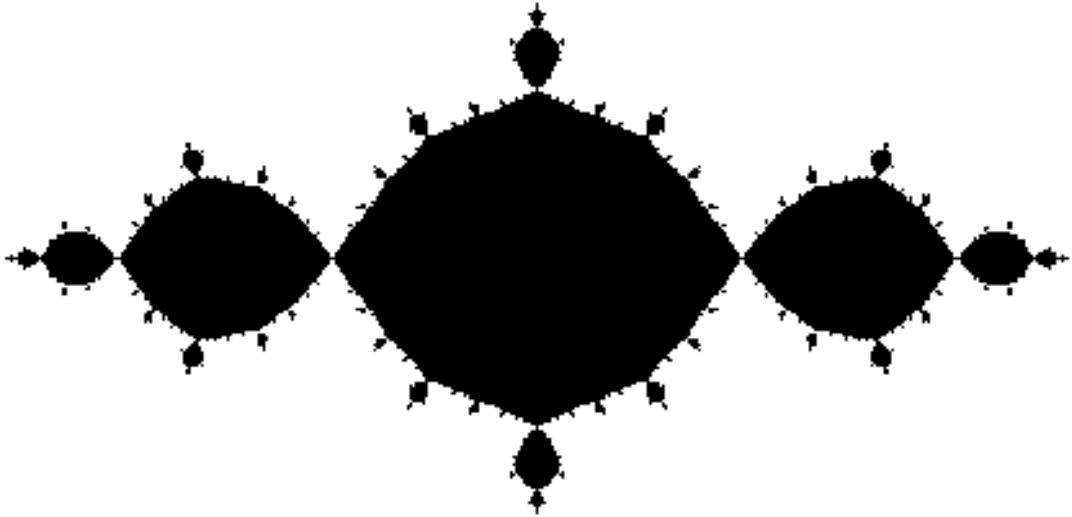,width=.9\hsize}}
\caption{The filled-in Julia set of $z \mapsto z^2-1$.}
 \label{c=-1.ps}  
\end{figure}

\begin{remark}  Concerning part 1 of Theorem A, we point out that in the 
paper \cite{TV}, Tangerman and Veerman showed that in the case of 
circle mappings with a non-flat singularity, the
parameter scaling and dynamical scaling agree for a 
large class of systems.   They have real methods comparing the  
dynamical derivatives and parameter derivatives along the critical value orbit.  Here, we use a complex technique for the unimodal scaling case since a direct derivative comparison appears to have extra difficulties.  
This is due to the changes in orientation, i.e., the folding which occurs for such 
maps, complicating the parameter derivative calculations.
\end{remark}

Figures \ref{Param puzzle} and \ref{c=-1.ps} illustrate item $2$ of the theorem.  
The reader is also encouraged to compare Theorem A with Theorem \ref{Lyu} of 
Lyubich on page \pageref{Lyu}. 

When dilating by the scaling factors given by the Fibonacci renormalization procedure, 
the computer pictures around the Fibonacci parameter also exhibit a hairy self-similarity. 
(Compare Figure \ref{Param puzzle} with Figure \ref{PM} on page \pageref{PM}.)  
Using the main construction of the proof of Theorem A, we demonstrate this hairiness.   
The appropriate scaling maps are denoted by $R_n$, and the Mandelbrot set by ${\bf M}$.

\begin{MT}{B:} {\bf (Hairiness for the Fibonacci parameter)}  

Given any disc $D(z,\epsilon)$ with center point $z$ and radius $\epsilon>0$, 
in the complex plane, there exists an $N$ such that for all $n > N$ we have that
\[  D(z,\epsilon) \cap R_n({\bf M}) \neq \emptyset.   \]
\end{MT}

The text is organized as follows.  In Section 1, we review some basic material of quadratic 
dynamics and the role of equipotentials and external rays.  In Section 2, we review the 
generalized renormalization procedure, where we define the principle nest in the 
dynamical plane as well as in the parameter plane.  In Section 3, we prove 
dynamical scaling and geometry results for the principal nest for parameter points 
which are Fibonacci renormalizable $n$-times.  These results and proofs are analogous to  
those given by Lyubich (\cite{Lg},  \cite{Lt}) for the Fibonacci point.  In Section 4, 
we construct a map of parameter space which allows us to compare it to the 
dynamical space and prove Theorem A, part 2.  In Section 5, we complete the 
proof of Theorem A.  Finally, in Section 6 we prove Theorem B.

\medskip

\noindent {\bf Acknowledgments:}  I would like to thank Misha Lyubich for our many insightful 
discussions and his continual encouragement.  
Thank you to Yair Minsky and John Milnor for making many
helpful suggestions for improving the exposition.  Also, thank you to Brian 
Yarrington for showing me how to create the pictures included 
in this paper, and to Jan Kiwi, Alfredo Poirier and 
Eduardo Prado for many interesting discussions.

\section{Introductory Material}
 
We outline some of the basics of complex dynamics of quadratic maps from the Riemann 
sphere to itself so that we may build the puzzle and parapuzzle pieces.  We will 
consider the normalized form, $f_c(z) = z^2 + c$ with parameter value $c \in \Bbb{C}$. 
The {\it basin of attraction} for infinity, $A(\infty)$, are all the points $z$ which 
converge to infinity under iteration.   The dynamics near infinity and the corresponding 
basin of attraction has been understood since B\"{o}ttcher (see \cite{Md}).  
Notationally we have that ${\Bbb D}_{r}$ is the disc centered at $0$ with radius $r$.

\begin{theorem}  
The map $f_c: A(\infty) \rightarrow A(\infty)$ is complex conjugate to the
map $w \mapsto w^2$ near infinity.  There exists a unique complex map  $\Phi_c$ defined 
on  $\widehat{\Bbb{C}} \setminus \bar{\Bbb{D}}_{r(c)}$, where $r(c)$ represents the 
smallest radius with the property that 
\[ f_c \circ \Phi_c(w) = \Phi_c(w^2), \] 
and normalized so that $\Phi_c(w) \sim w$ as $|w| \rightarrow \infty$.
\end {theorem} 

By Brolin \cite{Bro} the conjugacy map $\Phi_c$ satisfies 
\begin{equation} \lim_{n \rightarrow \infty} \log^+ (|f^n_c(z)|/2^n) = \log|\Phi^{-1}_c(z)|.
\label{Brolin} 
\end{equation}
In fact, the left hand side of equation (\ref{Brolin}) is defined for all 
$z \in A(\infty)$ and is the Green's function for the Julia set.

{\it Equipotential curves} are images of the circles with radii $r>r(c)$, 
centered at $0$ in $\widehat{\Bbb{C}} \setminus \bar{\Bbb D}_r$ under the map 
$\Phi_c$.  Actually, the moment that the above conjugacy breaks down is at the 
critical point $0$ if it is in $A(\infty)$.  In this case, if we try to extend 
the above conjugacy we see that the image of the circle with radius $r(c)$ 
passing through the critical point is no longer a disc but a ``figure eight".  
Despite the conjugacy difficulty, we may define equipotential curves passing 
through any point in $A(\infty)$ to be the level set from Brolin's formula.   
{\it External rays} are 
images of half open line segments emanating radially from ${\Bbb{D}}_{r(c)}$,  
i.e., $\Phi_c( re^{i \theta})$ with $r>r(c)$ and $\theta$ constant.  In fact, 
these are the gradient lines from Brolin's formula.  So again, we may extend 
these rays uniquely up to the boundary of $A(\infty)$ or up to where the ray meets the critical point or some preimage, i.e., the ``root" of a figure eight.

An external ray is referred to by its angle; for example the $\frac{1}{3}$-ray is 
the image of the ray with $\theta = \frac{1}{3}$.  A central question to ask is 
whether a ray extends continuously to the boundary of $A(\infty)$.  The following 
guarantees that some points (and their preimages) in the Julia set are such landing points. 

\begin{theorem} {\bf (Douady and Yoccoz, see \cite{Md})}         \label{per ray}
Suppose $z$ is a point in the Julia set which is periodic or preperiodic and the 
periodic multiplier is a root of unity or has modulus greater than $1$, then it is 
the landing point of some finite collection of rays.
\end{theorem}

One of the main objects of study in quadratic dynamics is the set of all parameters 
$c$ such that the conjugacy $\Phi_c$ is defined for the whole immediate basin of infinity. 

\begin{definition}
The {\it Mandelbrot set} {\bf M} consists of all values $c$ whose corresponding Julia 
set is connected.  
\end{definition}

The combinatorics of the Mandelbrot set have been extensively studied.  In \cite{DH3}, Douady and Hubbard present
many important results, some of which follow below.

\begin{theorem}{\bf (Douady and Hubbard, \cite{DH3})}  \label{dh3} \\                  
1.  The Mandelbrot set is connected. \\
2.  The unique Riemann map $\Phi_M: \widehat{\Bbb{C}} \setminus \bar{{\Bbb{D}}} 
\rightarrow  \widehat{\Bbb{C}} \setminus M$, with $\Phi_M(z) \sim z$ 
as $|z| \rightarrow \infty$, satisfies the following relation with the B\"{o}ttcher map:

\[ \Phi^{-1}_M(c) = \Phi^{-1}_c(c). \]
\label{phi}
\end{theorem}

With the Riemann map $\Phi_M$, we can define equipotential curves and external 
rays in the parameter plane analogous to the dynamical case above.  From the 
second result of Theorem \ref{phi}, it can be seen that the external rays and 
equipotentials passing through $c$ (the critical value) in the dynamical space 
are combinatorially the same external rays and equipotentials passing through 
$c$ in the parameter space.  Since the Yoccoz puzzle pieces have boundary which 
include rays that land, it is essential for the construction of the parapuzzle 
pieces that these same external rays land in parameter space.  Before stating 
such a theorem, we recall some types of parameter points.  {\it Misiurewicz 
points} are those parameter values $c$ such that the critical point of $f_c$ is 
pre-periodic.   A {\it parabolic point} is a parameter point in which the map 
$f_c$ has a periodic point with multiplier some root of unity.  For these points, 
their corresponding external rays land.

\begin{theorem}{\bf (Douady and Hubbard, \cite{DH3})}                   \label{misur}
If $c$ is a Misiurewicz point then it is the landing point of some finite collection
of external rays $R_M(\theta_i)$, where the $\theta_i$ represent the angle of the ray. 
In the dynamical plane, external rays of the same angle, $R_{f_c}(\theta_i)$, land at 
$c$ (the critical value of $f_c$).
\end{theorem}

\begin{theorem}{\bf (Douady and Hubbard, \cite{DH3})}                  
If $c$ is a parabolic point then it is the landing point of two external rays 
(except for $c=\frac{1}{4}$ which has one landing ray). In the dynamical plane for 
this $c$, these external rays land at the root point of the Fatou component 
containing the critical value.
\end{theorem}

Using rays and equipotentials in dynamical space, Yoccoz developed a kind of 
Markov partition, now called Yoccoz puzzle pieces, for non-renor\-mal\-izable 
(or at most finitely renormalizable) Julia sets with recurrent critical point 
and no neutral cycles \cite{H}.  Combining Theorems \ref{phi} and \ref{misur}, 
Yoccoz constructed the same (combinatorially) parapuzzle pieces for the parameter 
points of the non-renormalizable maps.

For our purposes we will now focus on the maps exhibiting initial behavior 
similar to the dynamics of the Fibonacci map.
The Fibonacci parameter value lies in what is called the $\frac{1}{2}$-wake.  
The $\frac{1}{2}$-wake is the connected set of all parameter values with 
boundary consisting of the $\frac{1}{3}$- and $\frac{2}{3}$-rays (which meet at 
a common parabolic point) and does not contain the main cardioid.   Dynamically, 
all such parameter points have a fixed point which is a landing point for the 
same angle rays, $\frac{1}{3}$ and $\frac{2}{3}$.  In fact, for all parameter 
points in the $\frac{1}{2}$-wake, the two fixed points are stable; we may follow them
 holomorphically in the parameter $c$.  We are now in a good position to 
review generalized renormalization in the $\frac{1}{2}$-wake.  We point out that this 
procedure, developed in  \cite{Lg}, is not restricted 
to the $\frac{1}{2}$-wake and the construction given below is readily generalized 
from the following description.

\section{A Review of Puzzles and Parapuzzles}

\noindent {\bf Initial Yoccoz Puzzle Pieces}

\begin{figure}[hbt]             
\centerline{\psfig{figure=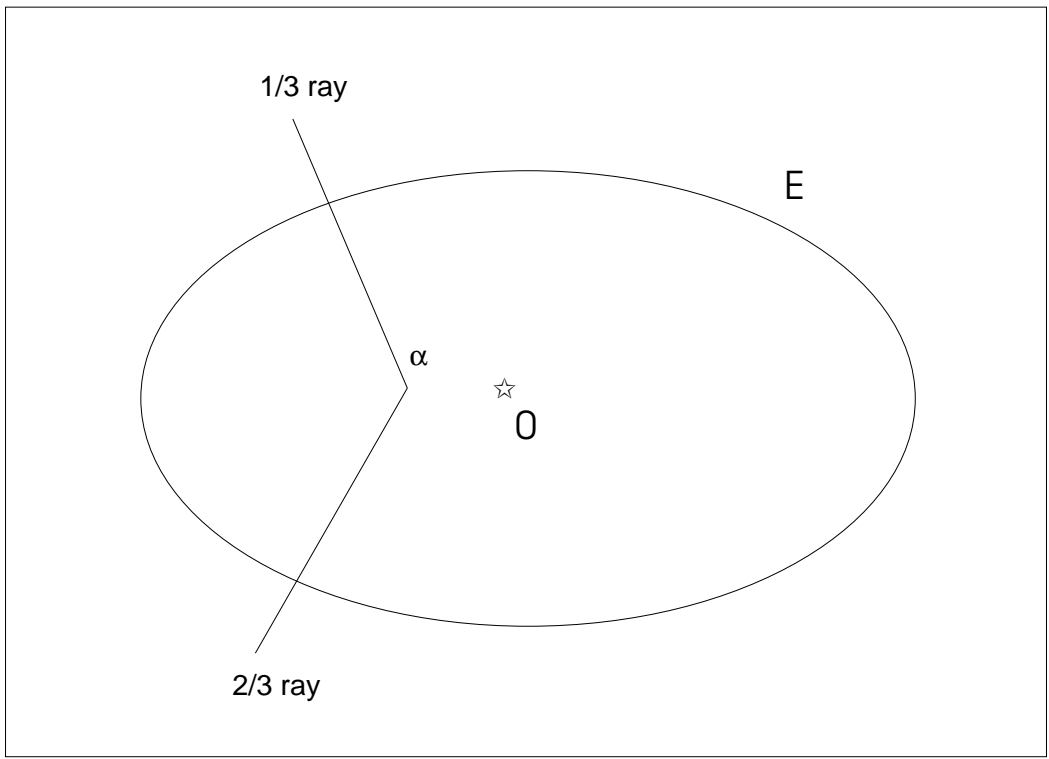,width=.9\hsize}}
\caption{Beginning generalized renormalization.}
\label{begin renorm}  
\end{figure}

We now review the Yoccoz puzzle piece construction essentially without proofs. 
(See \cite{H} or \cite{Lg} for more details.) For each parameter in the 
$\frac{1}{2}$-wake, we begin with the two fixed points commonly called the 
$\alpha$ and $\beta$ fixed points.  The $\beta$ fixed point is the landing point 
of the $0-ray$ (the only ray which maps to itself under one iterate).  The $\alpha$ 
point is the landing point of the $\frac{1}{3}$- and $\frac{2}{3}$-rays for all 
parameters in the $\frac{1}{2}$-wake.   By the B\"{o}ttcher map, it is easy to see 
that the $\frac{1}{3}$- and $\frac{2}{3}$-rays are permuted by iterates of $f_c$.  
The initial Yoccoz puzzle pieces are constructed as follows.  Fix an equipotential $E$.  
The top level Yoccoz puzzle pieces are the bounded connected sets in the plane with 
boundaries made up of parts of the equipotential $E$ and external rays. 
(See Figure \ref{begin renorm}.)  For the generalized renormalization procedure 
described below, the top level Yoccoz puzzle piece containing the critical point 
is labeled $V^0_0$.

\vspace{.2in}
\noindent {\bf The Principal Nest}

The generalized renormalization procedure for quadratic maps with recurrent 
critical point proceeds as follows.  For each parameter value $c$, iterate the 
critical point $0$ by the map $f_c$ until it first returns back to the set $V^0_0$.  
In fact, this will be two iterates. 
 Take the largest connected set around $0$, denoted $V^1_0$, such that 
$f^2(V^1_0) = V^0_0$.  Note that we suppress the parameter $c$ in this 
discussion, $V^n_0 = V^n_0(c)$.  This is the level $1$ central puzzle piece and 
we label the return map $f^2_c$ restricted to the domain $V^1_0$ by $g_{1}$ ($=g_{1,c}$).
  It is easy to see that $V^1_0 \subset V^0_0$ and that $g_1$ is a two-to-one branched 
cover.  The boundary of $V^1_0$ is made up of pieces of rays landing at points which 
are preimages of $\alpha$, as well as pieces of some fixed equipotential. 
 Now we proceed by induction.  Iterate the critical point until it first returns 
to $V^n_0$, say in $m$ iterates, and then take the largest connected set around $0$, 
denoted $V^{n+1}_0$. This gives $f^m_c(V^{n+1}_0) = V^n_0$.  Inductively we get a 
collection of nested connected sets $V^0_0 \supset V^1_0 \supset V^2_0 \supset V^3_0...$, 
and return maps $g_{n}(V^n_0) = V^{n-1}_0$.  Each of the $V^i_0$ has boundary equal to 
some collection of pieces of rays landing at preimages of $\alpha$ and pieces of some 
equipotential. Each $g_{i}$ is a two-to-one branched cover.  The collection of 
$V^i_0$ is called the {\it principal nest} of Yoccoz puzzle pieces around the critical point.

To define the principal nest of Yoccoz parapuzzle pieces in the parameter space it 
is easiest to view the above procedure around the critical value.  In this case, the 
principal nest is just the image of the principal nest for the critical point, namely 
$f_c(V^0_0) \supset f_c(V^1_0) \supset f_c(V^2_0) \supset f_c(V^3_0)...$.  Notice that 
again the puzzle pieces are connected and we have that the boundary of each puzzle piece 
to be some parts of a fixed equipotential and parts of some external rays landing at 
preimages of $\alpha$.  If we consider these same combinatorially equipotentials and 
external rays in the parameter space we get a nested collection of Yoccoz parapuzzle 
pieces.  By combinatorially the same we mean external rays with the same angle and 
equipotentials with the same values.

\begin{definition}
Given a parameter point $c$, the {\it parapuzzle piece of level $n$}, denoted by $P^n(c)$, 
is the set in parameter space whose boundary consists of the same (combinatorially) 
equipotentials and external rays as that of $f_c(V^n_0)$.
\end{definition}

We mention the essential properties about the sets $P^n$ used by Yoccoz. (The reader 
may wish to consult \cite{H} or \cite{GM}.)  The sets $P^n$ are topological discs.  
For all points $c$ in $P^n$,  the Yoccoz puzzle pieces of the principal nest (up to 
level $n$) are combinatorially the same.
This structural stability also applies to the off-critical pieces (up to level $n$) 
which are defined below. Hence, all parameter points in $P^n$ may be  renormalized in 
the same manner combinatorially up to level $n$.  We also point out that the set 
$P^n$ ($n>0$) intersects the Mandelbrot set only at Misiurewicz points.

\vspace{.2in}
\noindent {\bf  Off-critical Puzzle Pieces}

If a quadratic map is non-renormalizable, then at some level the principal nest is 
non-degenerate.  In other words, there is some $N$ such that for all $n \geq N$, 
$\mod (V^n_0,V^{n+1}_0)$ is non-zero. 
For these same $n \geq N$, we may iterate the critical point by the map $g_{n}$ some 
finite number of times until landing in $V^{n-1}_0 \setminus V^{n}_0$ (otherwise the 
map would be renormalizable).  Hence, to keep track of the critical orbit the generalized 
renormalization incorporates the following procedure.  Let us fix a level $n$.  For any 
point $x$ in the closure of the critical orbit contained in $V^{n-1}_0 \setminus V^{n}_0$, 
we iterate by $f_c$ until it first returns back to the $V^{n-1}_0$ puzzle piece. 
Denoting the number of iterates by $l$, we then take the largest connected neighborhood 
of $x$, say $X$, such that $f^{l}(X) = V^{n-1}_0$.  We only save those sets $X$, 
denoted $V^n_i$ ($i>0$), which intersect some point of the critical orbit.  We point 
out that the collection of $V^n_i$ are pairwise disjoint for $n>1$.  The return map 
$f^l(V^n_i)$ restricted to the set $V^n_i$ will still be denoted by $g_{n}$. 
The boundary of each $V^n_i$ must be a union of external rays landing at points which 
are some preimage of $\alpha$ and pieces of some equipotential.  Also, the return maps 
$g_{n}$ restricted to $V^n_i$ ($i \neq 0$) are univalent.  To review, for each level 
$n$ we have a collection of disjoint puzzle pieces $V^n_i$ and return maps, 

\begin{eqnarray*}
\bigcup_i V^n_i & \subset &  V^n_0, \\
g_{n} (V^n_i) &=& V^n_0. 
\end{eqnarray*}

\vspace{.2in}
\noindent {\bf  The Fibonacci Combinatorics}

Let us denote the Fibonacci sequence by $u(n)$, where $u(n)$ represents the $n$-th 
Fibonacci number.  The Fibonacci numbers are defined inductively: $u(0)=1, u(1)=1$ 
and $u(n)=u(n-1) + u(n-2)$.  The dynamical condition for $f_{c_{fib}}$ (recall $c_{fib}$ 
is real) is that for all Fibonacci numbers $u(n)$, we have 
$|f^{u(n)}(0)| < |f^{u(n-1)}(0)| < |f^i(0)|$, $u(n-1) < i < u(n)$.  
So the Fibonacci combinatorics require that the critical point return closest
to itself at the Fibonacci iterates. 

The generalized renormalization for the Fibonacci case is as follows.  (See  
\cite{LM} and \cite{Lt} for a more detailed account.  There is only one off-critical piece at every level, 
$V^n_1$.    The return map of $V^n_1$ to $V^{n-1}_0$ is actually just the restriction 
of the map $g_{n-1} : V^{n-1}_0 \rightarrow V^{n-2}_0$.  We point out that the map 
$g_{n-1}$ is the iterate $f^{u(n)}$ with restricted domain.  In short we have

\begin{eqnarray} 
g_{n} \: (\simeq g_{n-2} \circ g_{n-1}) & : & V^{n}_0 \rightarrow V^{n-1}_0  \: \:   (analytic \: double \: cover),     \nonumber  \\
g_{n-1} \:  & : & V^{n}_1 \rightarrow V^{n-1}_0 \: \: (univalent).  \nonumber 
\end{eqnarray}

\begin{figure}[hbt]        
\centerline{\psfig{figure=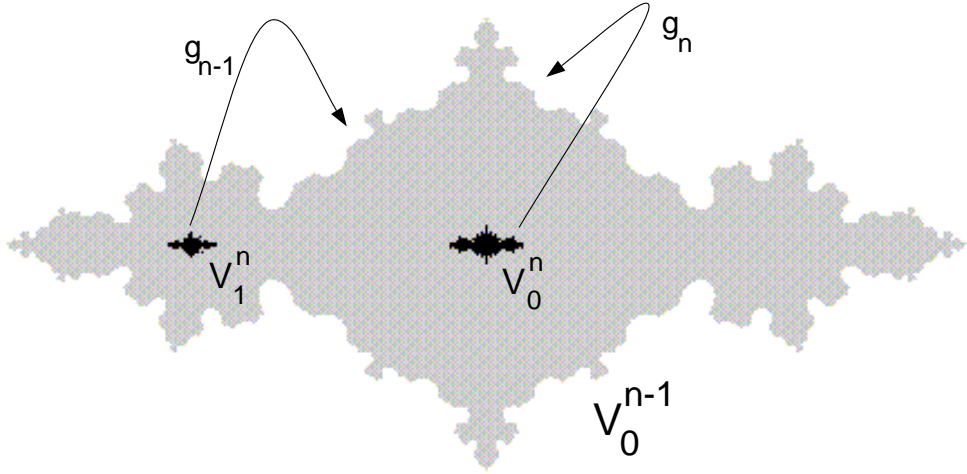,width=.9\hsize}}
\caption{Generalized renormalization: Fibonacci type return  with $n=7$.}
\label{Fib renorm}
\end{figure}

\begin{figure}[hbt]        
\centerline{\psfig{figure=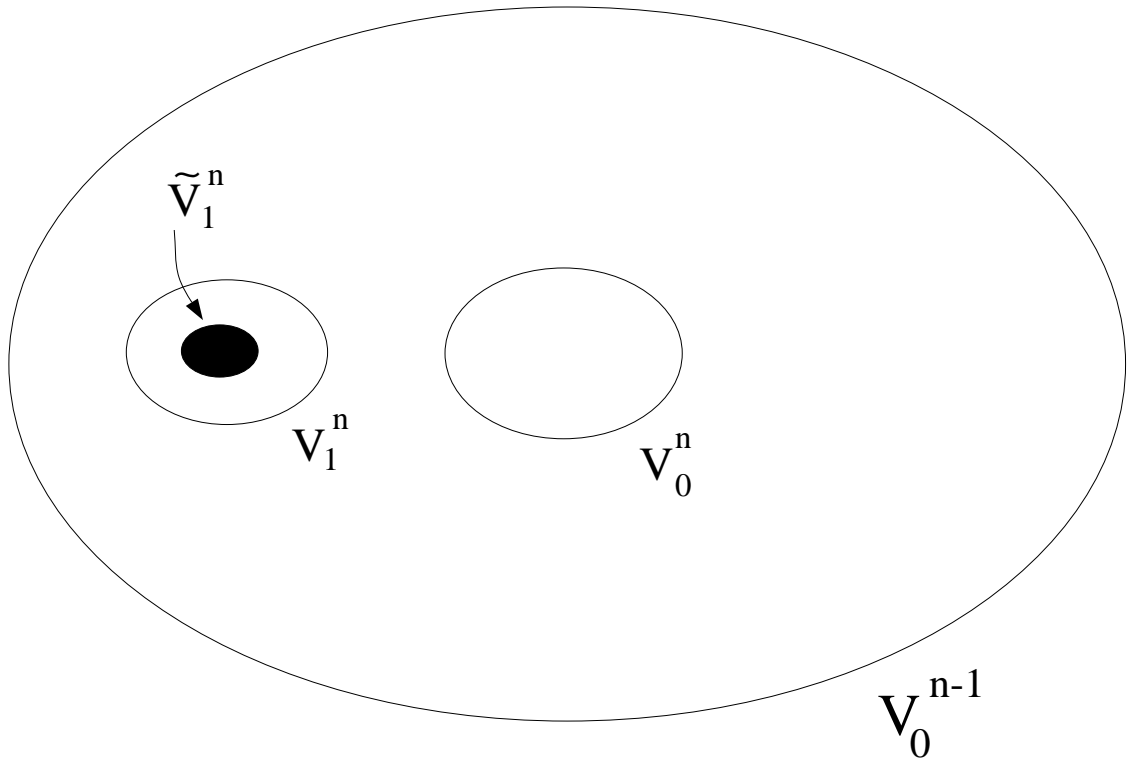,width=.9\hsize}}
\caption{Fibonacci puzzle piece nesting.}
\label{fib_ren}
\end{figure}

Finally, we define
the puzzle piece $\widetilde{V}^{n+1}_1$ to be the set $V^{n+1}_1$ which map
to the central puzzle piece of the next level down under $g_{n}$.  Namely, 
$\widetilde{V}^{n+1}_1$ is the set such that

\begin{eqnarray}
\widetilde{V}^{n+1}_1 &\subset &    V^{n+1}_1,     \nonumber \\
g_{n}(\widetilde{V}^{n+1}_1) &=& V^{n+1}_0. \nonumber 
\end{eqnarray}

\vspace{.2in}
\noindent {\bf Fibonacci Parapuzzle Pieces}

A Fibonacci parapuzzle piece, $P^n$, is defined as the set with the same combinatorial 
boundary as that of $f(V^n_0)$.  We also define an extra puzzle piece, $Q^n$.  In 
particular, $Q^n$ is a subset of $P^{n-1}$ and hence may be renormalized in the Fibonacci 
way $n-1$ times.  The boundary of the set $Q^n$ is combinatorially the same as 
$f(V^{n-1}_1)$.   Finally observe that $P^n \subset Q^n \subset P^{n-1}$.  
Properties for $P^n$ and $Q^n$ are given below. 

\begin{eqnarray}                 
c \in P^n  & \Longrightarrow & g_{n}(0) \in V^{n-1}_0   \nonumber \\
c \in Q^{n+1}  & \Longrightarrow & g_{n}(0) \in V^{n}_1 \nonumber \\
c \in P^{n+1}  & \Longrightarrow & g_{n}(0) \in \widetilde{V}^n_1  \nonumber 
\end{eqnarray}

We warn the reader that the parameter value $c$ has been suppressed as an index for the 
maps $g$ and puzzle pieces $V$.  Also, it is useful to use Figure \ref{Fib renorm} when 
tracing through the above properties of $Q^n, P^n$, and $P^{n+1}$, keeping in mind that 
$\widetilde{V}^n_1$, although too small for this picture, is contained in $ V^n_1$ 
(see Figure \ref{fib_ren}).

\vspace{.2in}
\noindent {\bf Lyubich's Motivating Result} 

The main motivating result of Lyubich is stated below.  We will give a brief review of 
the proofs and also point out that the proof of part $2$ of the theorem may be found 
in Lemma 4 of \cite{Lt}, while the proof of part $1$ is a direct consequence of part $2$ 
and the scaling results of \cite{Lg} (see pages 11-12 of this paper).  Finally, we point out that similar
scaling results were obtained on the real line in Lemma 5.4 of \cite{LM}.

\begin{theorem} {\bf (Lyubich)}           \label{Lyu}
The principal nest of central Yoccoz puzzle pieces for the Fibonacci map has the 
following properties.\\
1. The puzzle pieces scale down to the critical point in the following asymptotic manner:

\[ \lim_{n \rightarrow \infty}  \mod(V^{n-1}_{0}, V^{n}_{0}) \; / \; n  = \frac{1}{3} \ln 2.  \]  \\
2.  The rescaled puzzle pieces $V^n_0$ have asymptotic geometry equal to the filled-in 
Julia set of $ z \mapsto z^2 - 1$.  
\end{theorem}

The scaling factor in the theorem is exactly half that for the parameter scaling.  
This is because here the scaling is done around the critical point as opposed to the 
critical value.

\section{Beginning Geometry and Scaling}   \label{2?}

In studying the parameter space of complex dynamics, one first needs a strong command of 
the dynamics for all the parameter points involved.  Hence, before proceeding in the 
parameter space we shall first study the geometry of the central puzzle pieces, 
$V^n_0(c)$, for all $c \in Q^{n}$. 
 
Before stating a result similar to Theorem \ref{Lyu} for $c \in Q^n$, we indicate 
precisely how the rescaling of $V^n_0(c)$ is to be done.  For the parameter point 
$c_{fib}$, we dilate (about the critical point) the set $V^n_0(c_{fib})$ by a positive 
real constant so that the boundary of the rescaled $V^n_0(c_{fib})$ intersects the 
point $(1+\sqrt{5})/2$, the non-dividing fixed point for the map $z^2 - 1$.  If we dilate 
all central puzzle pieces $V^n_0(c_{fib})$ this way, we can then consider the rescaled 
return maps of $g_{n}$, denoted $G_{n}$, which map the dilated $V^n_0(c_{fib})$ to the 
dilated $V^{n-1}_0(c_{fib})$ as a two-to-one branched cover.  The map  $G_{n}$ restricted 
to the real line either has a minimum or maximum at the critical point.  To eliminate 
this orientation confusion, let us always rescale (so now possibly by a negative number) 
$V^n_0(c_{fib})$ so that the map $G_{n}$ always has a local minimum at the critical point.   

The point which maps to $(1+\sqrt{5})/2$ for $V^n_0(c_{fib})$ under this dilation we 
label $\beta_n$.  \label{bnc} Note that it must be some preimage of our original fixed 
point $\alpha$ and hence a landing point of one of the boundary rays.  (Puzzle pieces 
may only intersect a Julia set at preimages of $\alpha$.)  This point $\beta_n$ is 
parameter stable in that it may be continuously (actually holomorphically) followed for 
all $c \in Q^{n+1}$.  Hence we may write $\beta_n(c)$.  So for $c \in Q^{n+1}$, the 
rescaling procedure for $V^n_0(c)$ is to linearly scale (now possibly by a complex number) 
by taking $\beta_n(c)$ to $(1+\sqrt{5})/2$. 

The geometric lemma below gives asymptotic structure results for the central puzzle 
pieces $V^n_0(c)$ for all $c \in Q^n$.  In particular, the lemma indicates that as long 
as we can renormalize, the rescaled central puzzle pieces converge to the Julia set of 
$z^2-1$.  The notation $J\{z^2 - 1\}$ is used to indicate the Julia set for the map 
$z \mapsto z^2 - 1$.

Before proceeding we give a brief review of the Thurston Transformation (see \cite{DH2}) 
needed in the next lemma.  Consider the Riemann sphere punctured at 
$\infty, -1, 0$, and $\frac{1+\sqrt{5}}{2}$.  The map $\theta: z \mapsto z^2-1$ fixes 
$\infty$ and $\frac{1+\sqrt{5}}{2}$ while $-1$ and $0$ form a two cycle.  Consider 
any conformal structure $\nu$ on the Riemann sphere punctured at these points.  We 
can pull this conformal structure back by the map $\theta$.  This induces a map $T$ 
on the Teichm\"{u}ller space of the four punctured sphere.  A main result of this 
transformation $T$ is as follows.  

\begin{theorem} {\bf (Thurston)} Given any conformal structure $\nu$ we have that 
$T^n(\nu)$ converges at exponential rate to the standard structure in the 
Teich\-m\"{u}ller space. 
\end{theorem}

\begin{lemma} {\bf(Geometry of central puzzle pieces)}        \label{geom}
Given $\epsilon > 0$, there exists an $N>0$ such that for all $c \in Q^i$ where 
$i \geq N$ we have that the rescaled $\partial V^j_0(c)$ is $\epsilon$-close in the 
Hausdorff metric around $J\{z^2 - 1\}$, where $i \geq j+1 \geq N.$ 
\end{lemma}

\begin{proof}
First observe that the Julia set for $\Theta: z \mapsto z^2 - 1$ is hyperbolic.  
This means that given a small $\delta$-neighborhood of the Julia set there is some 
uniform contraction under preimages.  More precisely, there exists an integer $m$ 
and value $K>1$ such that for any point in the $\delta$-neighborhood of $J\{z^2 - 1\}$ 
we have 

\begin{equation}                    \label{contraction}
\max_{y \in \Theta^{-m}(x)} \dist (y, J\{z^2 - 1\}) < \frac{1}{K} \dist (x, J\{z^2 - 1\}).  
\end{equation}

Returning to our Fibonacci renormalization, it is a consequence of the main theorem of 
Lyubich's paper \cite{Lg} that the moduli of the nested central puzzle 
pieces, i.e., $\mod (V^n_0, V^{n-1}_0)$, grow at least at a linear rate independent 
of $c$.  Hence, independent of our parameters $c$ (although we must be able to renormalize 
in the Fibonacci sense), we have a definite growth in Koebe space for the map $g_{n}$. 
(This is somewhat misleading as the map $g_{n}$ is really a quadratic map composed with 
some univalent map.  Thus, when we say the map $g_{n}$ has a large Koebe space, we 
really mean that the univalent return map has a large Koebe space.)  The growing Koebe 
space implies that the rescaled maps $G_{n,c}$ have the following asymptotic behavior:

\begin{equation}       \label{oldG_n}   
G_{n,c}(z) = (z^2 + k(n,c))(1 + O(p^n)).     
\end{equation}
The bounded error term $O(p^n)$ comes from the Koebe space and hence, by the above 
discussion, is independent of $c$.  

We claim that $k(n,c) \rightarrow -1$ at an exponential rate in $n$, i.e., 
$|k(n,c) + 1|$ exponentially decays.  This result was shown to be true for 
$k(n,c_{fib})$ in Lemma 3 of \cite{Lt}.  We use this result as well as 
its method of proof to show our claim.  First, we review the method of proof used 
by Lyubich in the Fibonacci case.  This was to apply Thurston's transformation on 
the tuple $\infty$, $G_n(0)$, $0$, and $\frac{1+ \sqrt{5}}{2}$. 
 Pulling back this tuple by  $G_{n}$ results in a new 
tuple:  $\infty,$ the negative preimage of $G^{-1}_{n}(0)$, $0$, and 
$\frac{1+ \sqrt{5}}{2}$.  Next, two facts are used concerning the negative preimage of 
$G^{-1}_{n}(0)$.  The first is that it is bounded between $0$ and $- \frac{1+ \sqrt{5}}{2}$. 
(This is shown in \cite{LM}.) The second is that the puzzle piece $V^{n+1}_1$ is 
exponentially small compared to $V^{n}_0$ (a consequence of the main result of \cite{Lg}); 
therefore, after rescaling, the points $G^{-1}_{n}(0)$ and $G_{n+1}(0)$ are exponentially close. Hence, the tuple map 
\begin{equation}
(\infty, G_{n}(0), 0, \frac{1+ \sqrt{5}}{2}) \mapsto  ( \infty, G_{n+1}(0), 0, \frac{1+ \sqrt{5}}{2} ) 
\end{equation}
is exponentially close to the Thurston transformation since the pull-back by $G_{n}$ is 
exponentially close to a quadratic pull-back map (in the $C^1$ topology).   The Thurston 
transformation is strictly contracting; hence, the tuple must converge to its fixed point 
($\infty, -1, 0, \frac{1+ \sqrt{5}}{2}$).  Hence, we get $k(n,c_{fib}) \rightarrow -1$ at
 a uniformly exponential rate.  This concludes the summation of the Fibonacci case.

To prove a similar result for our parameter values $c$, let us choose some large level $n$ 
for which the $c_{fib}$ tuple is close to its fixed point tuple and such that the Koebe 
space for $g_{n}$ is large, i.e., $G_{n,c}$ is very close to a quadratic map.  Then we 
can find a small neighborhood around $c_{fib}$ in parameter space for which we still have 
a large Koebe space for $g_{n}$ (notationally $g_{n}$ is now $c$ dependent) and its 
respective tuple is also close to the fixed point tuple.  Then as long as the values $c$ 
in this neighborhood are Fibonacci renormalizable, we claim the $|k(n,c) + 1|$ 
exponentially decays.  We know that the Koebe space growth is at least linear and 
independent of the value $c$; hence by the strict contraction of the Thurston 
transformation we get our claim of convergence for $k(n,c)$.  Thus we may replace 
Equation (\ref{oldG_n}) with 

\begin{equation}        \label{G_n}
G_{n,c}(z) = (z^2 - 1)(1 + O(p^n)).  \nonumber   
\end{equation}

Returning to Equation (\ref{contraction}), we can state a similar contraction for the maps 
$G_{n,c}$.  In particular, in some small $\delta$-neighborhood of $J\{z^2-1\}$, we can find 
a value $k$ ($K>k>1$) and large positive integer $N_1$ so that for the same value $m$ as 
in Equation (\ref{contraction}) and for all $n>N_1$, we have

\begin{equation}                 \label{k}
\max_{y \in G^{-1}_{n+m-1,c} \circ G^{-1}_{n+m-2,c} \circ ... \circ G^{-1}_{n,c}(x) } \dist (y, J\{z^2 - 1\}) < \frac{1}{k} \dist (x, J\{z^2 - 1\}),  
\end{equation}
as long as $c$ is renormalizable in the Fibonacci sense, i.e., $n+m$ times.

From Equation (\ref{k}), we conclude the lemma.  We take the rescaled $V^n_0$ and note 
that it contains the critical point and critical value.  Hence we see that the topological 
annulus with boundaries $\partial{\Bbb D}_r$, $r$ large, and $\partial V^n_0$ under pull
 backs of $G_{n,c}$ must converge to the required set.  This concludes the lemma.
\end{proof}

The geometry of the puzzle pieces provides us with sufficient dynamical scaling results for 
the central puzzle pieces as well as for the off-critical puzzle pieces for $c \in Q^n$. 

\begin{lemma}                       \label{scale}
Given $\epsilon > 0$ there exists an $N$ so that for all $c \in Q^n$, $n>N$, we have the 
following asymptotics for the moduli growth of the principal nest 

\begin{equation}           \label{M-scale}
\left| \frac{\mod(V^{n}_{0}(c), V^{n+1}_{0}(c) )}{n} - \frac{1}{3} \ln 2 \right| < \epsilon.
\end{equation}

\end{lemma}

\begin{proof}
Notationally we will suppress the dependence of the parameter $c$.   By Lyubich 
(\cite{Lg}, page 12), the moduli growth from  $\mod (V^{n-1}_{0}, V^{n}_{0})$ to 
$\mod (V^{n}_{0}, V^{n+1}_{0})$ approaches 
$\frac{1}{2}(-\capp_{\infty}(J\{z^2-1\}) - \capp_{0}(J\{z^2-1\}))$. 
(See Appendix for the definition of capacity.) The proof of the growth relies only 
on the geometry of the puzzle pieces.  The map $g_n$ takes the annulus 
$V^{n}_{0} \setminus V^{n+1}_{0}$ as a two-to-one cover onto the annulus 
$V^{n-1}_{0} \setminus \widetilde {V}^n_1$.  Hence we have  the equality

\begin{equation}                     \label{V, tilde V}
\mod (V^{n}_{0}, V^{n+1}_{0}) = \frac{1}{2} \mod (V^{n-1}_{0}, \widetilde {V}^n_1). 
\end{equation}

Using the Gr\"{o}tzsch inequality on the right hand modulus term, we have 

\begin{equation}             \label{V's}
\mod (V^{n-1}_{0}, \widetilde {V}^n_1) =  \mod (V^{n-1}_{0}, V^{n}_{1})  + \mod(V^{n}_{1}, \widetilde {V}^n_1) + a(V^{n-1}_{0}, \widetilde {V}^n_1), 
\end{equation}
where the function $a(V^{n-1}_{0}, \widetilde {V}^n_1)$ represents the Gr\"{o}tzsch error. 
By applying the map $g_{n-1}$ to $V^n_1$ we see that $\mod (V^{n}_{1}, \widetilde {V}^n_1)$ is equal to $\mod (V^{n-1}_{0}, V^{n}_{0})$.  The term 
$\mod (V^{n-1}_{0}, V^{n}_{1})$ converges to $\mod (V^{n-2}_{0}, V^{n-1}_{0})$.  
This is easily seen by applying the map $g_{n-1}$ which is a two-to-one branched cover
with the critical point image being pinched away from $V^{n-1}_0$ as 
$n \rightarrow \infty$.  Finally, the Gr\"{o}tzsch error $a$ depends only on the geometry of 
$V^n_0$ because of the linear increase in modulus  between both $V^{n-1}_0$ and 
$V^{n+1}_0$ and hence is approaching $-\capp_{ \infty}(V^n_0) - \capp_{0}(V^n_0)$ (see Lemma 
\ref{caps} in the Appendix).  But this is approaching $-\capp_{\infty}(J\{z^2-1\}) - 
\capp_{0}(J\{z^2-1\})$ by Lemma \ref{geom} and the fact that the capacity function preserves 
convergence in the Hausdorff metric (see Lemma \ref{Schiffer}).  Finally, 
$-\capp_{\infty}(J\{z^2-1\}) - \capp_{0}(J\{z^2-1\})$ is shown to be equal to $\ln 2$ in the 
Appendix.  Using the notation $m_n = \mod (V^{n-1}_{0}(c), V^{n}_{0}(c))$, we may 
rewrite Equation (\ref{V's}) as
\[ m_{n+1} = \frac{1}{2}m_n + \frac{1}{2}m_{n-1} + \frac{1}{2}a + o(1), \]
where $a = \ln 2$. The asymptotics of this equation give the desired result. 
\end{proof}

\section{The Parameter Map}

\noindent {\bf Dynamical puzzle piece rescaling}

Now that we have a handle on the geometry of the central puzzle pieces for values $c$ in our 
parapuzzle, let us consider rescaling the $V^n_0$ in a slightly different manner.  For 
each $c \in Q^n$ dilate $V^{n}_0$ so that the point $g^{-1}_{n}(0)$ maps to $-1$.  
Notice this is just an exponentially small perturbation of our previous rescaling since 
there we had $G^{-1}_n(0)$ approaching $-1$ uniformly in $n$ for all $c \in Q^n$.  Hence 
Lemmas \ref{geom} and \ref{scale} still hold for this new rescaling.  Let us denote this 
new rescaling map by $r_{n,c}$.  Therefore, fixing $c \in Q^n$, the map $r_{n,c}$ is 
the complex linear map $x \mapsto (1 / g^{-1}_{n}(0)) \cdot x$. 

\begin{lemma}                           \label{univ}
The rescaling map $r_{n,c}$ is analytic in $c$.  In other words,  $g^{-1}_{n}(0)$ is 
analytic in $c \in Q^n$.
\end{lemma}

\begin{proof}
The roots of any polynomial vary analytically without branching provided no two collide.  
We claim the root in question does not collide with any other.  But for all $c \in Q^n$ we 
have that the piece $V^{n+1}_1(c)$ can be followed univalently in $c$.  Hence, we have our 
claim. 
\end{proof} 

We remind the reader that the map $g_{n,c}$ is just a polynomial in $c$.  Let us define 
the analytic parameter map which allows us to compare the dynamical space and the 
parameter space. 

\smallskip 

\noindent {\bf The Parameter Map:}  The map $M_n(c)$ is defined as the map 
$c \mapsto r_{n,c} \cdot g_{n+1,c}(0)$ with domain $c \in Q^n.$ 

Since the map $r_{n,c}$ is just a dilation for fixed $c$, we see that if $M_n(c)=0$ 
then this parameter value must be superstable. This superstable parameter value, denoted 
$c_n$, is the unique point which is Fibonacci renormalizable $n$ times, and for the 
renormalized return map, the critical point returns precisely back to itself, 
i.e., $g_{n}(0)=0$.  Equivalently, this is the superstable parameter whose  critical point 
has closest returns at the Fibonacci iterates until the $n+1$ Fibonacci iterate when 
it returns to itself, $f^{u(n+1)}_{c_n}(0)=0$.

\begin{lemma} {\bf (Univalence of the parameter map.)}    \label{ML}
For sufficiently large $n$, there exists a topological disc $S^n$ such that 
$P^n \subset S^n \subset Q^n$, the map $M_n(c)$ is univalent in $S^n$, and 
$\mod(S^n, P^n)$ grows linearly in $n$.
\end{lemma}

The proof of Lemma \ref{ML} is technical so we give an outline for the reader's 
convenience.  We first show that the winding number is exactly $1$ around the image 
$-1$ for the domain $P^n$.  This will be a consequence of analysis of a finite number 
of Misiurewicz points along the boundary of $P^n$.  Using Lemma \ref{geom} we will 
locate the positions (up to some small error) these selected Misiurewicz points must 
map to under $M_n(c)$.  Then we prove that the image of the segments in 
$\partial P^n$ between these Misiurewicz points is small, where ``between" is defined 
by the combinatorial order of their rays and equipotentials.  Hence, the 
$c \in \partial P^n$ have to follow the combinatorial order of the points of 
$J\{z^2-1\}$ without much error.  Since we wind around $-1$ only once when traveling 
around $J\{z^2-1\}$ the only way we could have more than one 
preimage of $-1$ for the map $M_n(c)$ would be for one of these segments of $\partial P^n$ to 
stretch a ``large" distance and go around the point $-1$ a second time.  But this 
cannot happen if the segments follow the order of $J\{z^2-1\}$ without much error.  
Finally, we show that this degree one property extends to some increasingly large 
image around $-1$ in Lemma \ref{big space}.

\begin{proof}
We will again use the map $\Theta(z)=z^2-1$.  Let $b_0 = \frac{1+\sqrt{5}}{2}$ be the 
non-dividing fixed point for the Julia set of $\Theta$.  The landing ray for this 
point is the $0$-ray.  Taking a collection of pre-images of $b_0$ under the map 
$\Theta$ we may order them by the angle of the ray that lands at each point. 
(Note that there is only one angle for each point.)  The notation for this combinatorial 
order of preimages will be $b_0, ... b_i, b_{i+1}, ..., b_0$.  

Since the point $b_0$ is in the Julia set of $\Theta$, the set of all preimages of 
$b_0$ is dense in the Julia set.  Given that this Julia set is locally connected we 
have the following density property of the preimages of $b_0$:  given any 
$\epsilon > 0$, we can find an $l$ so that the collection of preimages 
$\Theta^{-l}(b_0)$ is such that the Julia set between any two successive 
points (in combinatorial order) is compactly contained in an $\epsilon$-ball.  
In other words, for this set $\Theta^{-l}(b_0)$, given any $b_i$ and $b_{i+1}$, 
the combinatorial section of the Julia set of $\Theta$ between these two points is 
compactly contained in an $\epsilon$-ball.

For each $c \in Q^n$ we define an analogous set of points $b_{i,n}(c)$ along the 
boundary of the rescaled puzzle pieces $V^n_0(c)$.  First let us return to our old 
way of rescaling $V^n_0(c)$, taking the point $\beta_n(c)$ to $1+\sqrt{5}/2$ 
(see page \pageref{bnc}).  For our value $l$ above we take a set of points to be 
preimages of $(1+\sqrt{5})/2$ under the map $G_{n-1} \circ ... \circ 
G_{n-l}$ for each $c$.  These points are on the boundary of the rescaled $V^n_0(c)$ 
and in particular are endpoints of some of the landing rays which make up some of the 
boundary of the rescaled $V^n_0$.  In particular, we may label and order this set 
of preimages $b_{i,n}(c)$ by the angles of their landing rays.  Hence we may also 
refer to a piece of the rescaled boundary of $V^n_0(c)$ as a piece of the boundary 
that is combinatorially between two successive $b_{i,n}(c)$.

We claim that for $n$ large enough we have that for all $c \in Q^n$ these 
combinatorial pieces of $r_{n,c}(V^n_0)$, say from $b_{i,n}$ to $b_{i+1,n}$, is 
in the exact same $\epsilon$-ball as their $b_{i}$ to $b_{i+1}$ piece counterpart.  
For this claim we first want $b_{i,n}(c) \rightarrow b_i$ as $n \rightarrow \infty$.  
But this is true (for this rescaling) by the proof of Lemma \ref{geom} since the 
rescaled maps $G_n$ converge to $\Theta$ exponentially. 

Now that we have a nice control of where the Misiurewicz points of $\partial P^n$ are 
landing, we focus on the boundary segments of $P^n$ between them.  Note that by 
Theorem \ref{dh3} of Douady and Hubbard, we have a good combinatorial description of 
$\partial P^n$ in terms of rays and equipotentials.  Combinatorially the image of 
these boundary segments under the map $M_n$ will be in the appropriate boundary 
segments of the dynamical puzzle pieces.  Therefore, we focus on 
controlling the combinatorial segments between the $b_{i,n}$ along the central 
puzzle pieces in dynamical space.  With precise information on where these 
combinatorial segments are in dynamical space we make conclusions on the image 
of $\partial P^n$.  

Now we prove that the combinatorial piece between $b_{i,n}(c)$ and $b_{i+1,n}(c)$ 
converges to the combinatorial piece from $b_i$ to $b_{i+1}$ in the Hausdorff metric.  
Let us take a small neighborhood around $c_{fib}$ such that the rescaled $V^n_0$ 
are in some small neighborhood around $J\{z^2-1\}$.  For all $c$ in this neighborhood, 
take the combinatorial piece $b_{i,n}(c)$ to $b_{i+1,n}(c)$ such that the distance 
(in the Hausdorff metric) is greatest from $b_i$ to $b_{i+1}$.   
Suppose this distance is $\delta$, then after $m$ preimages (the value $m$ being the 
same as in Lemma \ref{geom}, see Equations (\ref{contraction}) and (\ref{k})), 
the distances between these preimages is less than $\delta/\lambda^m$, where $\lambda>1$ 
and is independent of the parameter.  Finally, notice that for the $b_i$ segments, 
any preimages of a combinatorial segment must be contained in another $b_i$ segment 
(the Markov property).  Hence, we actually get convergence at an exponential rate.

To review, the points $c \in \partial P^n$ under the map $M_n(c)$ must traverse 
around the point $-1$ with each appropriate Misiurewicz point landing very near 
$b_i$ since for all $c$, $b_{i,n}(c) \rightarrow b_i$.  But each combinatorial piece 
is also very near the combinatorial piece for the Julia set of $\Theta$ and the 
Julia set has winding number $1$ around the point $-1$ which completes the winding 
number argument for this rescaling.  Now if we rescale by $r_{n,c}$ instead of the 
old way (they are exponentially close) the same result holds. This completes the proof of 
the univalence of the map at least in some small image containing $-1$.  The lemma 
below will complete the proof of this lemma.
\end{proof}

\begin{lemma}                                   \label{big space}
For all sufficiently large $n$, there exists $R(n) \rightarrow \infty$ as 
$n \rightarrow \infty$ such that the map $c \mapsto M_n(c)$ is univalent 
onto the disc $D(-1,R(n))$.
\end{lemma}

\begin{proof}
The image of any point $c \in \partial Q^n$ under the map $M_n(c)$ is contained in the set 
$r_{n,c}(V^{n-1}_0(c))$.  But the boundary of $r_{n,c}(V^{n-1}_0(c))$ under the rescaling of 
$r_{n,c}$ is very far from $r_{n,c}(V^{n}_0(c))$ by the modulus growth proven in 
Lemma \ref{scale} (see Appendix, Proposition \ref{far} and reference).  Let $R$ equal 
the minimum distance from the image of $\partial Q^n$ to the origin.  Note 
that $Q^n \setminus P^n$ cannot contain the point $-1$ in its image under $M_n(c)$ 
since the closest these points can map to $-1$ is when they map into a small neighborhood 
of $J\{z^2-1\}$.  Since we showed in the proof above that the winding number around $-1$ 
for $M_n(\partial P^n)$ is one, we must have the same result for $M_n(\partial Q^n)$ 
since $-1$ can have no new preimages in this domain $Q^n \setminus P^n$.  Hence, the 
winding number is one for all points in the disc of radius $R$.  Thus, the map $M_n$ 
must be univalent in some domain with image (at least) the disc centered at $0$ and 
radius $R$.  Taking the preimage of this disc will define the desired set in parameter 
space, $S^n$.  The result follows and hence does Lemma \ref{ML}.
\end{proof}

Lemma \ref{ML} also allows us to give the geometric result of Theorem A.
As 
$n$ increases we have an increasingly large Koebe space around the image of $\partial P^n$.  
Since the image of $\partial P^n$ under the map $M_n(c)$ must asymptotically approach 
that of $J\{z^2-1\}$, the parapuzzle pieces must also
asymptotically approach this same geometry by application of the Koebe Theorem.  
Hence, by Lemma \ref{ML}, we get part 2 of Theorem 
A.

\section{Parapuzzle Scaling Bounds}

\indent To understand the scaling in parameter space, we focus on the image of the 
parapuzzle pieces $P^n$ and $P^{n+1}$ under the parameter map $M_n$.  Since $M_n$ is 
nearly a linear map for the domain $P^n$,
we are in a good position to prove the scaling results of the Main Theorem A.

\begin{theorem} {\bf (Theorem A,  part 1.)}
The principal nest of Yoccoz parapuzzle pieces $P^n$ for the Fibonacci point $c_{fib}$ 
scale down in the following asymptotic manner:
\[  \lim_{n \rightarrow \infty}  \mod (P^{n}, P^{n+1}) \;  / \;  n  =  \frac{2}{3} \ln 2.  \]
\end{theorem}

\begin{proof}
We begin by defining two bounding discs for the $J\{z^2-1\}$.  Take as a center the 
point $0$ and fix a radius $T$ so that the disc $D(0,T)$ compactly contains $J\{z^2-1\}$.  
Also take a radius $t$ so that the disc $D(0,t)$ is strictly contained in the immediate 
basin of $0$ for $J\{z^2-1\}$. (See Figure \ref{z^2-1 centered}.)  This gives

\begin{equation}               \label{J center} 
J\{z^2-1\} \subset D(0, T ) \setminus D(0, t).
\end{equation}

Let us calculate the scaling properties of the image of $\partial P^{n+1}$ under the 
same map $M_n$.  Again we will have that the image of $\partial P^{n+1}$ ``looks'' 
like $J\{z^2-1\}$ although at a much smaller scale.  We remind the reader that $M_n$ 
maps the point $c_{n+1}$ to $-1$.  Now we claim that the point $c_{n+1}$ acts as 
the ``center" of $\partial P^{n+1}$ in the following sense:

\begin{equation}               \label{p center}
M_n( \partial P^{n+1} ) \subset D(-1, \frac{M^{\prime}_n}{M^{\prime}_{n+1}} T ) \setminus D(-1, \frac{M^{\prime}_n}{M^{\prime}_{n+1}} t),
\end{equation}
where $M^{\prime}_n$ represents the derivative of $M_n$ at the point $c_{n}$.  To prove 
the claim we note that $M_{n+1}( \partial P^{n+1} ) \subset  D(0, T ) \setminus D(0, t)$. 
Pulling this image back by the univalent map $M_n \circ M^{-1}_{n+1}$ and noting that 
this map has increasing Koebe space for our domain proves this claim.

\begin{figure}[hbt]          
\centerline{\psfig{figure=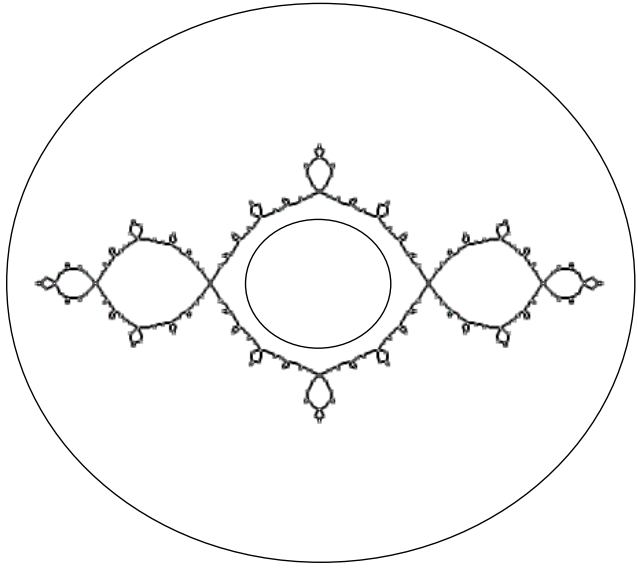,width=3.3in}}
\caption{The centering property for the Julia set of $z^2-1$.}
\label{z^2-1 centered}
\end{figure}

Now let us observe what is happening dynamically for all $c \in P^{n+1}$.  We have 
that  $r_{n,c}( \widetilde{V}^n_1)(c)$ is also centered around $-1$ by the construction 
of $r_{n,c}$.  Hence we have a result similar to that in expression (\ref{J center}), 
although perhaps with different radii.  Most importantly, however, the different radii 
must preserve the same centering ratio seen in expression (\ref{p center}), 
i.e., $\frac{T}{t}$. 

To compare the centerings of the dynamical and parameter sets above, we focus on 
the Fibonacci point $c_{fib}$.  We have that the point $M_n(c_{fib})$ is contained 
in the topological annulus of expression (\ref{p center}).   But this image must 
also be contained in the centering annulus of $r_{n,c_{fib}}(\partial \widetilde{V}^n_1)$ 
in the dynamical space.  Geometrically the point $M_n(c_{fib})$ is to the sets 
$r_{n,c_{fib}}(\partial \widetilde{V}^n_1)$ and $M_n( \partial P^{n+1} )$ as the 
$-1$ point is to the Julia set of $z \rightarrow z^2-1$ up to exponentially small error.  
Hence we have the following equivalent centerings

\begin{equation}               \label{centering2a}
r_{n,c_{fib}}(\partial \widetilde{V}^n_1)(c_{fib}) \subset D( -1, \frac{M ^ { \prime } _ n }{M^{\prime}_{n+1}} T) \setminus D(-1, \frac{M^{\prime}_n}{M^{\prime}_{n+1}} t),
\end{equation}

\begin{equation}               \label{centering2b}
M_n( \partial P^{n+1} ) \subset 
D(-1, \frac{M ^ { \prime } _ n }{M^{\prime}_{n+1}} T) \setminus D(-1, \frac{M^{\prime}_n}{M^{\prime}_{n+1}} t).
\end{equation}
Let us rewrite the scaling estimate of Equation (\ref{V, tilde V}) from Lemma \ref{scale},  

\[ \lim_{n \rightarrow \infty}  \mod \left( V^{n-1}_0(c_{fib}), \widetilde{V}^n_1(c_{fib}) \right) \; / \; n  = \frac{2}{3} \ln 2. \]
Since the modulus function is preserved under rescalings, we apply $r_{n,c_{fib}}$ to get

\begin{equation}          \label{scale cor}
\lim_{n \rightarrow \infty}   \mod \left( r_{n,c_{fib}}(V^{n-1}_0(c_{fib})), r_{n,c_{fib}}( \widetilde{V}^n_1(c_{fib}) ) \right) \; / \; n  = \frac{2}{3} \ln 2.
\end{equation}  
Expressions (\ref{centering2a}) and (\ref{centering2b}) and Equation (\ref{scale cor}) 
combined with Lemma \ref{ML} give  

\begin{equation}              
\lim_{n \rightarrow \infty} \mod \left( M_n(P^n),  M_n(P^{n+1}) \right) \; / \;  n = \frac{2}{3}  \ln 2,
\end{equation}
which proves Theorem A, part 1, and hence completes the proof of this
theorem.
\end{proof}

\section{Hairiness at the Fibonacci Parameter}    

\indent Let us define the Mandelbrot dilation for the Fibonacci point given by the 
renormalization.   We wish to dilate the Mandelbrot set, {\bf M}, about the Fibonacci 
parameter point by taking the approximating superstable parameter points $c_{n}$ to 
some fixed value for each $n$.  Of course, we have been doing a similar kind of dilation 
in the previous section so we will take advantage of this work and rescale in 
the following more well-defined manner.   

\medskip

\noindent {\bf Mandelbrot rescaling:}         \label{f-hairy}
Let $R_n$ be the linear map acting on the parameter plane which takes $c_{fib}$ to 
$-1$ and $c_{n}$ to $0$.  Notice that this is nearly the same map as our parameter map 
$M_n$.  The maps $M_n$ have an increasing Koebe space, take $c_{n}$ to $0$, and 
asymptotically takes $c_{fib}$ to $-1$. 

\medskip 
 
The proof of hairiness will be a consequence of the geometry of the external rays
which make up pieces of the boundary of the principal nest puzzle pieces, $V^n_0(c)$.  
Before proving this theorem, we first give a combinatorial description of how these 
rays lie in the dynamical space for the Fibonacci parameter.

We remind the reader that $\beta_{n,0}$ is on the boundary of $V^n_0$ and is the 
landing point of two external rays.  We label the union of these two rays of 
$\beta_{n,0}$ as $\gamma(\beta_{n,0})$.  The curve $\gamma(\beta_{n,0})$ divides 
the complex plane into two regions.  We label the region which does not contain 
the piece puzzle $V^n_0$ as $\Gamma(\beta_{n,0})$.  

We also define similar objects $\Gamma(x)$ and $\gamma(x)$ for the other Julia set 
points $x$ on the boundary of $V^n_0$.  To start, we have the symmetric point 
$\beta_{n,1}$ of $\beta_{n,0}$, and note $g_{n}(\beta_{n,1}) = g_{n}(\beta_{n,0})$.  
We can exhaust all other Julia set points on the boundary of $V^n_0$, denoting them as 
$\beta_{n,i}$ where $g_{n-i+1} \circ g_{n-i-1} \circ... \circ g_{n} (\beta_{n,i}) = 
\beta_{n-i,0}$ for $2 \leq i \leq n$.  Of course this representation is not unique in 
the variable $i$ but we will not need to distinguish between these various 
$\beta_{n,i}$ points.   For each of the $\beta_{n,i}$ points we can define 
$\gamma(\beta_{n,i})$ as the union of the two external rays which land there.  
Similarly we define the $\Gamma(\beta_{n,i})$ region as we did for $\beta_{n,0}$.  
In particular, $\Gamma(\beta_{n,i})$ has boundary $\gamma(\beta_{n,i})$ and does 
not contain $V^n_0$.   

The combinatorial properties for the $\gamma$ and $\Gamma$ sets are easy to determine 
for the Fibonacci parameter.  First we have that $|\beta_{n,0}| < |\beta_{n-1,0}|$ 
where the absolute values are necessary since the $\beta$'s change orientation 
(see page \pageref{bnc}).  If the $\beta_{n,0}$ and $\beta_{n-1,0}$ have the same 
sign then  $\Gamma(\beta_{n,0}) \supset \Gamma_{n-1}(\beta_{n-1,0})$, otherwise we 
replace $\beta_{n,0}$ with its symmetric point to achieve this inclusion.  By application 
of pull-backs of $g_n$ it is easy to see that

\begin{equation}    \label{beta}
\bigcup_{i} \Gamma(\beta_{n,i}) \supset \bigcup_{i} \Gamma(\beta_{n-1,i}). 
\end{equation}

Since this is just a combinatorial property depending on the first $n$ Fibonacci 
renormalizations, this property holds as we vary our parameter $c$ in $Q^n$.  As 
a direct consequence of expression (\ref{beta}), we conclude that

\begin{equation}    \label{Jc}
J_c \cap (V^{n-1}_0(c) \setminus V^n_0(c)) \subset \bigcup_{i} \Gamma_n(\beta_{n,i}(c)). 
\end{equation}

\begin{figure}[htp]          \centerline{\psfig{figure=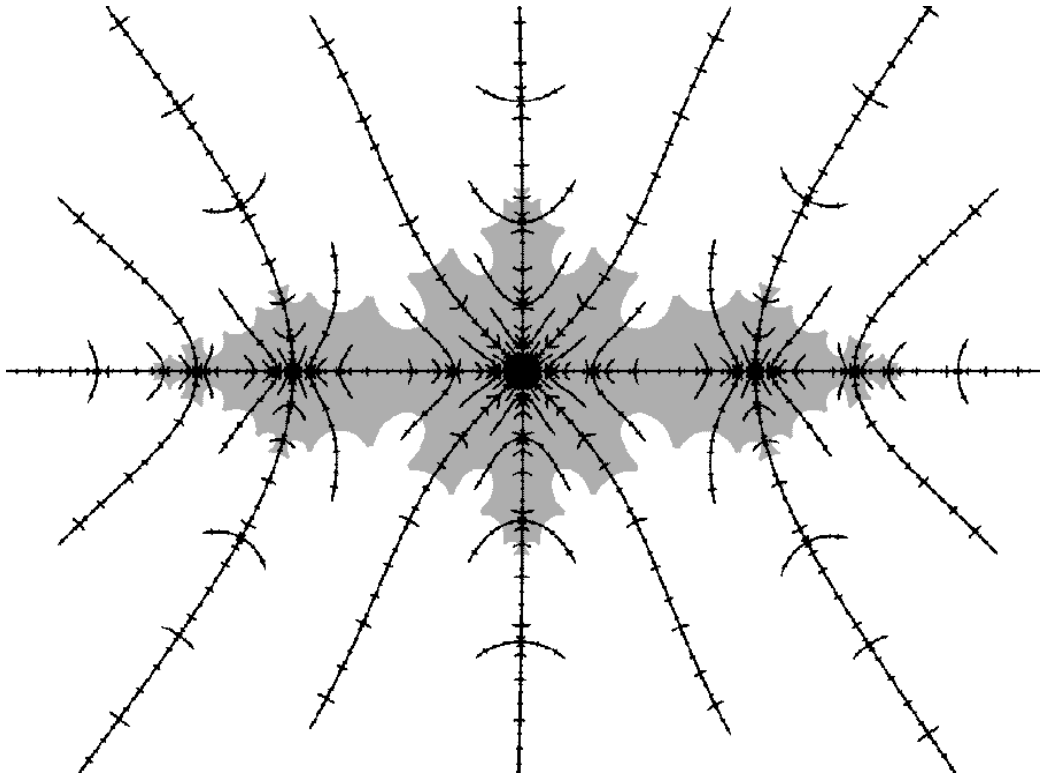,height=.45\vsize}}
\caption{Parapuzzle piece $P^6$ with Mandelbrot set.}
\label{PM}
\end{figure}

\begin{figure}[htp]          \centerline{\psfig{figure=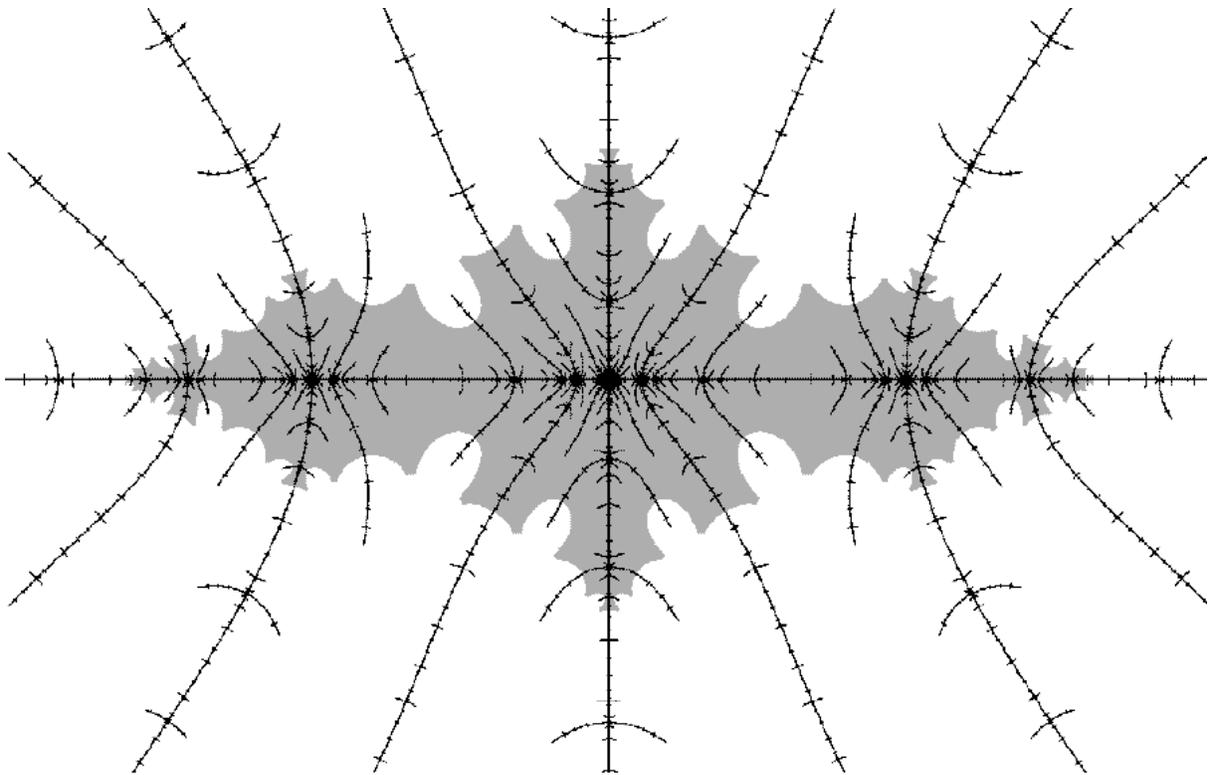,height=.45\vsize}}
\caption{Dynamical central puzzle piece $V^6_0$ for the Fibonacci Julia set.}
\label{PJ}
\end{figure}

\begin{figure}[hbtp]        
\centerline{\psfig{figure=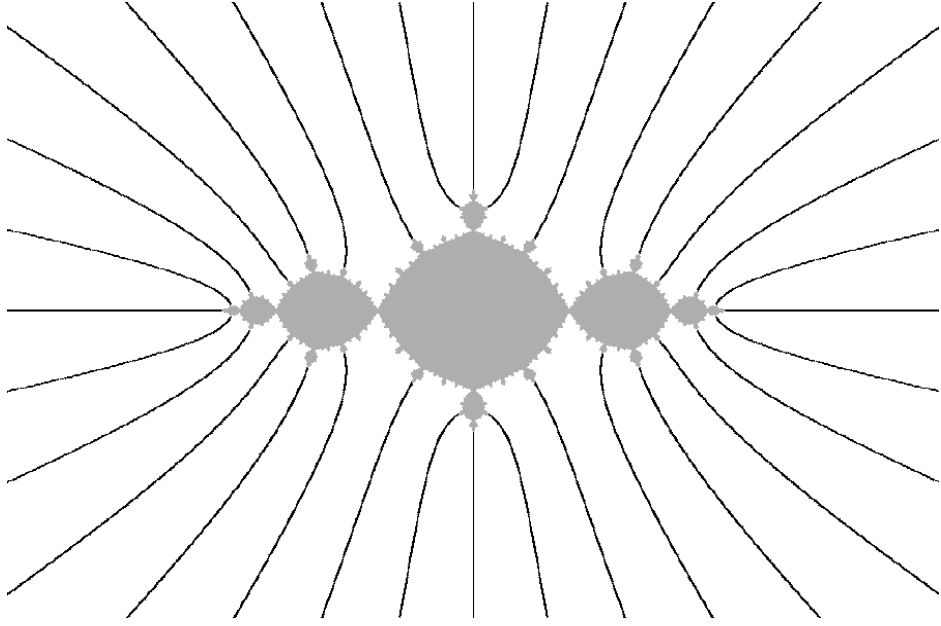,width=.75\hsize}}
\caption{The $0$-ray and some of its preimages for the Julia set of $z^2-1$.}
\label{0-rays-pic}
\end{figure}

By the dynamical scaling results we know that if we rescale the left side of 
expression (\ref{Jc}) by $r_{n,c}$ then $\partial V^{n-1}_0$ tends to infinity 
while $V^n_0(c)$ stays bounded (see Appendix, Proposition \ref{far}).  Hence for
 connected Julia sets the appropriate connected pieces must ``squeeze through" 
the $\Gamma$ regions in $V^{n-1}_0 \setminus V^n_0$.    We will be able to conclude 
the hairiness theorem by application of our map $M_n$ and by the geometry of the 
$\Gamma$ regions, i.e., the controlled ``hairiness" of $J_c$.  We show that the 
rescaled $\Gamma$ regions, i.e.,  $r_{c,n}(\Gamma_n(\beta_{n,0}(c)))$ are converging 
to the $0$-ray of $J\{z^2-1\}$.  (Compare Figures \ref{PM} and \ref{PJ} with 
\ref{0-rays-pic}.)  Also we show that $r_{c,n}(\Gamma_{n+1}(\beta_{n+1,0}(c)))$ 
converges to the inner $0$-ray of the Fatou component containing $0$ for $J\{z^2-1\}$.

\begin{lemma}                 \label{rotate}
For $c \in P^n$ the linear rescaling maps $r_{n,c}$ and $r_{n,c_{fib}}$ have 
asymptotically the same argument, 
$| \text{arg} (r_{n,c}) - \text{arg} (r_{n,c_{fib}}) | \rightarrow 0$ modulo $\pi$.
\end{lemma}

\begin{proof}
The return maps $g_{n}$ are asymptotically $z^2-1$ post-composed and pre-composed by 
a linear dilation.  When our return maps have a large Koebe space we see that the 
rescaling argument difference (as in the Lemma) converges to a constant modulo $\pi$.  
For the Fibonacci parameter case we are always rescaling by a real value so the 
difference is $0$ modulo $\pi$.  Since we are scaling down to the Fibonacci parameter 
we get the desired result.
\end{proof}

\begin{lemma}                 \label{gamma}
For discs $D(0,\rho)$ in the plane, there exists an $N(\rho)>0$ so that for all 
$n > N$ the curves $r_{n,c}(\Gamma(\beta_{n,0}(c)))$ converge to the $0$-ray of the 
Julia set of $z \mapsto z^2-1$ in the Hausdorff metric in $D(0,\rho)$.  Also, the 
curves $r_{n,c}(\Gamma(\beta_{n+1,0}(c)))$ converge to the inner $0$-ray  of the 
Fatou component containing $0$ for the Julia set of $z \mapsto z^2-1$. 
\end{lemma}

\begin{proof}

By Lemma \ref{rotate} the rescaling maps $r_{n,c}$ converge to a real dilation.   
Hence there is a decreasing amount of ``rotation" in the return map $g_{n,c}$.  In 
particular, the return maps $g_{n}$ are close to $z^2-1$ post-composed and pre-composed 
with a real rescaling in the $C^1$ topology.  Let us focus on the curves 
$r_{n,c}(\Gamma(\beta_{n,0}(c)))$.  Since we know that the pull-backs are 
essentially $z^2-1$, the curves should converge as stated in the theorem.  
However, there are two difficulties.  First, our $G_n$ pull-backs are not defined 
in all of ${\Bbb C}$ and second, $z^2-1$ is contracting under preimages.  Hence, 
we check that after pulling back our curves $r_{n,c}(\Gamma(\beta_{n,0}(c)))$ by 
$G_n$ that their extensions (i.e., the rescaled pull-back of the whole curve by the 
appropriate $f$ iterate) have some a priori bounds.  

Let us take the set $\gamma_{n} \cap V^{n-2}_0$ and pull-back by $g_{n} \circ g_{n-1}$.  
Taking the appropriate branches we get 
$(\gamma(\beta_{n+1,0}) \cup \gamma(\beta_{n+1,1}) \cup \gamma(\beta_{n+1,2} )) \cap V^{n}_0$.
In particular, the endpoints of $\gamma(\beta_{n+1,0})$ lie on the boundary of $V^{n}_0$.  
Hence their extension is determined by property (\ref{beta}) (the geometry of the rays 
of the previous level).  In particular, we have that 
$\gamma(\beta_{n+1,0}) \cap ( V^{n-1}_0 \setminus V^{n}_0 )$ is combinatorially 
between $\Gamma(\beta_{n,0})$ and $\Gamma(\beta_{n,2})$.  The piece of 
$\gamma(\beta_{n+1,0})$ contained in $V^{n}_0$ is controlled by the nearly 
$z^2-1$ pull-backs (the maps $g_{n-1} \circ g_{n-2}$ after rescaling) of again 
lesser level rays as constructed above.  So let us assume the sets $\gamma(\beta_{j,0})$, 
$\gamma(\beta_{j,1})$ and $\gamma(\beta_{j,2})$, $j = \{n, n-1\}$ nicely lie in 
the appropriate half-planes, where nice means that $r_{c,n} (\gamma(\beta_{j,0}))$ is 
in the right-half plane, $r_{c,n} (\gamma(\beta_{j,1}))$ in the left-half plane, and 
$r_{c,n} (\gamma(\beta_{j,2}))$ in the upper-half plane.  Then by the above argument 
we have that the collection $r_{n+1,c} (\gamma(\beta_{j,0}))$, 
$r_{n+1,c} (\gamma(\beta_{j,1}))$ and $r_{n+1,c} (\gamma(\beta_{j,2}))$ with 
$j = n+1$ is also nice in that they lie in the appropriate half-planes.  
This completes the induction step. 
 
The initial step comes from the fact that the geometry is nice in the Fibonacci case.  
More precisely we have that $r_{c,n} (\gamma(\beta_{n,0}))$ is contained in the 
right-half plane just by symmetry.  If we pull-back as above we see that 
$r_{c,n} (\gamma(\beta_{n,0}))$ must be contained in the right-half plane.  
Hence we may perturb this set-up in a small parameter neighborhood to start 
the induction process.

Because the return maps $G_n$ uniformly (in parameter $c$) approach $z^2-1$, we may 
use the a priori bounds and the coordinates from the B\"{o}ttcher map of $z^2-1$ 
to conclude that the rescaled rays $r_{c,n}(\Gamma(\beta_{n,0})$ must uniformly 
approach the $0$-ray of $z^2-1$ in compact sets.  Finally, viewing this same pull-back 
argument inside of $r_{n,c}(V^n_0)$ for the curves $r_{n,c}(\Gamma(\beta_{n+1,0}(c)))$ 
yield convergence to the inner $0$-ray and completes the lemma. 

\end{proof}

We are now in a good position to prove hairiness in an arbitrary disc 
$D(z,\epsilon) \subset {\Bbb {C}}$.   We point out that if $z$ is in  
$J\{z^2-1\}$, the theorem holds by Lemma \ref{ML}.   In this lemma we showed that 
the Misiurewicz points on the boundary of $P^n$ under our rescaling map, $M_n$, 
converge to the preimages of the $\beta$ fixed point of $z^2-1$.   Note that the 
preimages of the $\beta$ fixed point are dense in $J\{z^2-1\}$.  Given that the 
map $M_n(c)$ is an exponentially small perturbation of $R_n(c)$ we must have hairiness 
for neighborhoods of such $z$ and this claim is proven.   In fact, the above argument 
shows that it suffices to show that images $M_n({\bf M})$ satisfy the Theorem B.

\medskip

\noindent {\bf Proof of hairiness: }

\begin{proof} 
We first focus on the structure of $J_c$ for parameters $c$ in $Q^n$.  By 
Lemma \ref{gamma}, we have that $r_{c,n}(\Gamma(\beta_{n,0}(c)))$ converges 
to the $0$-ray of $J\{z^2-1\}$ in bounded regions.  Hence, for $c \in {\bf M} \cap P^n$ 
we must have that its Julia set in this region, 
i.e., $r_{c,n}(\Gamma(\beta_{n,0}(c))) \cap J_c$, also converges to the $0$-ray  
(compare property (\ref{Jc})).  Now the image of $M_n( {\bf M} \cap Q^n )$ must map 
into the set $\cup_c r_{c,n}(\Gamma(\beta_{n,0}(c))) \cap J_c)$.  Also, this domain 
contains the Misiurewicz point, say $c^{\prime}$, which lands at the rescaled $\beta$ point 
$r_{c^{\prime},n}(\beta_{n-1,0})$.  But $ r_{c^{\prime},n}(V^{n-1}_0) $ is growing at an exponential 
rate while $r_{c,n}(\beta_{n,0})$, the ``other" end of this image, converges to the 
$\beta$ fixed point of $z^2-1$ for all $c \in Q^n$.  Note we must have a Misiurewicz 
point landing near this $\beta$ point as well.  Because the Mandelbrot set is connected 
we get that a piece of the image $M_n( {\bf M} \cap Q^n )$ converges to the $0$-ray 
of $J\{z^2-1\}$.  Similarly we have convergence of 
$\cup_c r_{c,n}(\Gamma(\beta_{n+1,0}(c)) \cap J_c)$ to the inner $0$-ray of $J\{z^2-1\}$. 
Hence pieces of $M_n( {\bf M} \cap Q^n )$ also have convergence to this inner $0$-ray.

So given an arbitrary disc $D(z,\epsilon)$, we iterate it forward by $z^2-1$ until it 
intersects the $0$-ray or inner $0$-ray of the Julia set of $z \mapsto z^2-1$.   By 
the above we have that this image will eventually intersect all Julia sets of 
$P^n \cap {\bf M}$.  Pulling back by our almost $z \mapsto z^2-1$ maps shows that all 
Julia sets $P^n \cap {\bf M}$ must eventually intersect $D(z,\epsilon)$.  Applying our 
parameter map and arguing as above yields hairiness.
\end{proof}

\appendix
\section{Geometry of sets in the plane } 

\noindent {\bf Topological discs in the plane. }

We define capacity for sets in the plane and reference the perturbation result used 
in this paper.  We point out that there are many equivalent definitions of capacity, 
many of which may be found in the book of Ahlfors (chapter 2, \cite{A}).  We give one 
such definition.  Take a topological disc, $U$ in the plane with boundary 
$\partial {U} = \Gamma$.  Fix a point $z \in U$.  Let $\cal R$ be the Riemann map 
of the unit disc onto $U$ with $\cal R$$(0) = z$.

\begin{definition}
The {\it capacity} of $U$ (or $\Gamma$) with respect to the point $z$ is  
\[  cap_z(U) = \ln {\cal R}^{\prime}(0).   \]
\end{definition}

We can calculate the capacities needed for this paper.  For $cap_{\infty}(J\{z^2-1\})$ 
we proceed as follows.  Using the B\"{o}ttcher map and Brolin's formula, we see that 
the dynamics for the attracting basin is conjugate to the complement of the unit disc 
under the $z \mapsto z^2$ map.  The conjugacy is in fact the Riemann mapping which 
has derivative precisely $1$ at infinity (in the appropriate coordinate system). 
Hence, $cap_{\infty}(J\{z^2-1\}) = \ln 1 = 0$.  

Similarly, we may calculate $cap_{0}(J\{z^2-1\})$.  (Note we must only consider the 
connected component containing $0$ for the capacity definition.)  The dynamics around 
the critical point $0$ is $z \mapsto 2z^2$ (two iterates of $z \mapsto z^2-1$).  Again 
we can conjugate the immediate basin of attraction for the critical point to 
$z \mapsto z^2$  with domain the unit disc (again by the B\"{o}ttcher map).  Comparing 
the two maps,  $z \mapsto 2z^2$ and  $z \mapsto z^2$, we see that the conjugacy 
(Riemann map) must have derivative equal to $1/2$. Hence 
$cap_{0}(J\{z^2-1\}) = \ln \frac{1}{2}$.

To state a perturbative result of capacity, consider all topological disc 
boundaries $\Gamma$ in the plane with the Hausdorff metric $d_H$.  The following 
result says that if we fix a point $z$ bounded away from some $\Gamma_{\infty}$, 
then exponential convergence to this curve in the Hausdorff metric yields 
exponential convergence in their capacities.  The result is due to Schiffer and 
may be found in his paper \cite{SCH} or the book of Ahlfors \cite{A} pages 98-99.

\begin{theorem} {\bf (Schiffer)}      \label{Schiffer}
Given a sequence of disc boundaries $\Gamma_i$ with convergence at an exponentially 
decreasing rate to some $\Gamma_{\infty}$, 
$ d_H(\Gamma_n,\Gamma_{\infty}) = O(p^n) $, and a point $z$ bounded away from 
$\Gamma_{\infty}$, then $cap_z(\Gamma_n) = cap_a(\Gamma_{\infty}) + O(p^n)$
\end{theorem}

\vspace{.1in}

\noindent {\bf Topological annuli in the plane. }

We take two topological open discs $U_1$ and $U_2$ in the plane such that $U_2$ 
is compactly contained in  $U_1$.  Then we may form the annulus 
$A = U_1 \setminus \bar{U}_2$.  Every such annulus can be mapped 
(a {\it canonical map}) univalently to an annulus $ \{ z | \; 0 < r_1 < |z| < r_2 \}$.   
Although an annulus can be mapped to many different such annuli, there does exist a 
conformal invariant, namely the ratio of the radii $ r_2 / r_1 $.  There are many 
equivalent definitions for the modulus of an annulus, one of which is given here.

\begin{definition}
The {\it modulus} of an annulus $A$, $mod A$, is the conformal invariant 
$\log \frac{r_2}{r_1}$ resulting from a canonical map. 
\end{definition}

\begin{MT}  {\bf (Koebe: Analytic version)} 
Take any two topological discs $U_1$, $U_2$ with ${U}_2 \subset \subset U_1$, and 
a univalent map $g$ with domain $U_1$.  Then independently of the map $g$, there exists `
a constant $K$ such that 

\[  \frac{|g'(x)|}{|g'(y)|}  < K      \]
for $x,y \in U_2$.  Also $K = 1 + O(\exp(-mod(U_2 \setminus U_1)))$ as 
$mod(U_2 \setminus U_1) \rightarrow \infty$.
\end{MT}

\begin{MT} {\bf (Gr\"{o}tzsch Inequality)}  Given three strictly nested topological `
discs, $U_3 \subset U_2 \subset U_1$ in $\Bbb{C}$,

\[ mod(U_1 \setminus {U}_2) +   mod(U_2 \setminus {U}_3) \leq   mod(U_1 \setminus {U}_3). \]

\end{MT}

Now suppose we take a sequence of $U_3(i)$, containing $0$ and converging to $0$, and 
a sequence $U_1(i)$ with boundary converging to infinity.  The set $U_2$ will remain fixed.  
Also, suppose  $U_3(i)  \subset U_2 \subset U_1(i) \subset \Bbb{C}$, then the 
equipotentials for the topological annuli $U_1(i) \setminus U_3(i)$ in compact regions 
of ${\Bbb C} \setminus 0$ converge to circles centered at $0$.  One consequence is the 
following proposition.

\begin{proposition}    \label{caps}
Given $U_1(i)$, $U_2$, and $U_3(i)$, the deficit in the Gr\"{o}tzsch Inequality converges 
to $cap_0(U_2) + cap_{\infty}(U_2). $
\end{proposition}

Finally, we mention one extremal situation (see \cite{LV}, Chapter 2 for actual estimates). 
Suppose we take a topological annulus $A \in {\Bbb{C}}$ with inner boundary 
$\Gamma_1$ and outer boundary $\Gamma_2$.

\begin{proposition}  \label{far}
If $A$ is normalized so that the diameter of $\Gamma_1$ is equal to $1$ then 
$\dist(\Gamma_1, \Gamma_2) \rightarrow \infty$ as $ \mod A  \rightarrow \infty$.
\end{proposition}


 \bibliography{Lthesis}
 \bibliographystyle{amsalpha}

\end{document}